
\magnification=\magstep1 
\baselineskip=14pt

\def\operator#1#2{\def#1{\mathop{\kern0pt\fam0#2}\nolimits}} 
\operator\ini{in}
\operator\Gin{Gin} 
\operator\GGin{{ Gins }} 
\operator\GL{GL}
\operator\Ln{Ln}
\operator\Support{Support}
\operator\mm{ {\bf m}}
\operator\revlex{RL}
\operator\depth{depth}
\operator\Lex{Lex}
\operator\Tor{Tor}
\operator\Ker{Ker}
\operator\CoKer{CoKer}
\operator\Image{Image}
\operator\chara{char}
\operator\rank{rank}
\operator\Borel{Borel}
\operator\projdim{projdim}
\operator\reg{reg}

\def\NN{ {\bf N} }

\def\RR{ {\bf R} }
 
\def\endProof{ \hbox{}\nobreak  
\vrule width 1.6mm height 1.6mm depth 0mm \par
\goodbreak \medskip}

\def\nt{\noindent}

\font\big=cmbx10 scaled\magstephalf
\ \vskip 3cm

\centerline {\big  Koszul homology and extremal properties of Gin and  Lex }

\vskip 0.5cm 
\centerline {Aldo Conca} \centerline{ Dipartimento di Matematica, Universit\'a di Genova}
\centerline{Via Dodecaneso 35, I-16146 Genova, Italia}
\centerline{conca@dima.unige.it}

\vskip1cm

 \centerline {1. \bf Introduction} 
\medskip

 Let $K$ be a field of characteristic $0$ and $I$ be a homogeneous ideal of the polynomial ring
$R=K[x_1,\dots,x_n]$. We   denote by $\beta_{ij}(R/I)$ and by $\beta_i(R/I)$  the graded Betti numbers and the total
Betti numbers  of
$R/I$, that is: 

$$\beta_{ij}(R/I)=\dim_K  \Tor^R_i(R/I,K)_j \quad \hbox{ and  }
\quad \beta_{i}(R/I)=\dim_K 
\Tor^R_i(R/I,K)$$  where the subscript $j$ on the right  of a graded module denotes, throughout the paper, the degree
$j$ component of that module.
   
There are two monomial ideals   canonically attached   to
$I$: the generic initial ideal $\Gin(I)$ with respect to the revlex order and the lex-segment ideal
$\Lex(I)$.  They  play a fundamental role in the investigation of many algebraic,  homological, combinatorial and
geometric  properties of the ideal 
$I$ itself.  By definition, the generic initial ideal $\Gin(I)$ is the initial ideal  of $I$  with respect to the revlex
order  after performing a generic change of coordinates.  The ideal $\Lex(I)$ is defined as follows. For every vector 
space $V$ of forms of degree, say, $d$ one defines $\Lex(V)$ to be the vector space generated by the largest, in the
lexicographic order,  $\dim V$  monomials of degree $d$.  For a homogeneous ideal $I$ one sets $\Lex(I)=\oplus_{d}
\Lex(I_d)$. By the very definition,
$\Lex(I)$ is simply a graded vector space but Macaulay's theorem on Hilbert functions, see for instance [V, Sect.1],  says
that $\Lex(I)$ is indeed an ideal.  By construction,  it is clear that $\Lex(I)$ only depends on the Hilbert function of
$I$.  The graded Betti numbers  of $I, \Gin(I)$ and $\Lex(I)$    satisfy the following inequalities: 
\medskip

\nt {\bf Theorem 1.1} \ \  {\sl
 (a) $\beta_{ij}(R/I)\leq  \beta_{ij}(R/\Gin(I))$  for all $i,j$, \par
\nt (b) $\beta_{ij}(R/I) \leq \beta_{ij}(R/\Lex(I))$  for all $i,j$. }
\medskip 

The first  inequality can be proved by a standard deformation argument and holds, more generally,   when 
$\Gin(I)$ is replaced by any initial ideal of $I$, generic or not. The second
   is an important theorem proved in characteristic $0$ by Bigatti [B] and Hullett [H]  (independently)  and  extended
later by Pardue [P2] to positive characteristic.

Ideals having the same Betti numbers as their generic initial ideal are characterized by the  following result of 
Aramova, Herzog  and  Hibi [AHH, Theorem 1.1]: 
\medskip

\nt {\bf Theorem 1.2} \ \  {\sl
 The following conditions are equivalent:
\par
\nt (a) $\beta_{ij}(R/I)=  \beta_{ij}(R/\Gin(I))$ for all $i,j$, \par
\nt (b) $\beta_{1j}(R/I) =  \beta_{1j}(R/\Gin(I))$ for all $j$, \par
\nt (c) $\beta_{1}(R/I) =  \beta_{1}(R/\Gin(I))$, \par
\nt (d) $I$ is componentwise linear.  }
\medskip 

Conditions (b) and (c)  are not explicitly stated in [AHH] but it follows from the proof that they are indeed  
equivalent to the others. Similarly,  for inequality  (b) of Theorem 1.1,    Herzog and Hibi [HH, Corollary 1.4] proved: 
\medskip 

\nt {\bf Theorem 1.3}   \ \  {\sl
 The following conditions are equivalent:
\par
\nt (a) $\beta_{ij}(R/I) =  \beta_{ij}(R/\Lex(I))$ for all $i,j$, \par
\nt (b) $\beta_{1j}(R/I)  =   \beta_{1j}(R/\Lex(I))$ for all $j$, \par
\nt (c) $\beta_{1}(R/I)  =   \beta_{1}(R/\Lex(I))$, \par
\nt (d) $I$ is a Gotzmann ideal. } 
\medskip 

Recall that a homogeneous ideal $I$ is said to be {\it componentwise linear} if for all 
$k\in \NN$ the ideal $I_{<k>}$ generated by the elements of degree $k$ of $I$ has a linear resolution.  Also,
$I$ is said to be a {\it Gotzmann}  ideal  if for every
$k\in
\NN$ the space
$I_k$ of forms of degree $k$ in $I$   has the smallest possible span in the next degree according to the Macaulay
inequality [V, Thm3.1], i.e.,  $\dim_K R_1I_k=\dim_K R_1\Lex(I_d)$.   

Betti numbers can be computed via Koszul homology. Let $K(y, R/I)$ be the graded Koszul complex with values in
$R/I$ where $y=y_1,\dots,y_n$ is a system of generators of the maximal homogeneous ideal of $R$  and let
$H_i(y,R/I)$ be its homology. Then 
$$\beta_{ij}(R/I)=\dim _K H_i(y,R/I)_j. $$ 

More generally,  we can consider Koszul homology with respect to a sequence of generic linear forms. We define: 
\medskip

\nt {\bf Definition 1.4}   \  { For $p=1,\dots,n$, let  $K(p,R/I)$ be the graded Koszul complex with    respect to a
sequence of $p$ generic linear forms 
$y_1,\dots,y_p$. Denote by
$H_i(p,R/I)$ the $i$-th homology  of $K(p,R/I)$.  These are graded modules. We define the {\it Koszul-Betti numbers} of
$R/I$  by setting 
$$\beta_{ijp}(R/I)=\dim_K H_i(p,R/I)_j.$$  We extend the definition   to $p=0$ by setting
$H_0(0,R/I)=R/I$ and also
$H_i(p,R/I)=0$ whenever $i>p$.   }
\medskip

We  show in this paper that  Theorems 1.1, 1.2, 1.3 hold more generally for  Koszul-Betti numbers.   This is  done in
Section 4 while Sections 2 and  3 contain some preliminary results. 

In Section 5 we investigate the properties of the set $\GGin(I)$ of all the generic initial ideals of an ideal
$I$.   We   show that the Koszul-Betti numbers of the revlex gin, $\Gin(I)$,  are less than or equal to than those
of any other gin of $I$. On the other hand, we exhibit  examples of ideals for which  in $\GGin(I)$  there is  no ideal
whose Betti numbers are
greater than or equal to    those of any other element of $\GGin(I)$.  To construct these examples we introduce the
notion of almost Borel-fixed ideals. The main feature of  an almost Borel-fixed ideal $I$ is that the set
$\GGin(I)$ can be described in a quite simple way.  Then the Eliahou-Kervaire formula for the Betti numbers of
Borel-fixed ideals guides us rapidly to the construction of the examples mentioned above. In Section 6 we list some open
questions related  to the results of Sections 4 ad 5. 
\bigskip\bigskip\bigskip

\centerline{2. \bf Preliminaries} 
\medskip

In this section  we introduce some definitions, notation and some preliminary facts.  For generalities  on term orders,
Gr\"obner bases, initial ideals and lex-segments  we refer the reader to [E,BH,G,GS,KR,S,V].  In the following we will
consider only term orders
$\tau$ such that 
$$x_1>_\tau x_2 >_\tau \dots >_\tau x_n.$$

By definition, the $\beta_{ijn}(R/I)$'s  are the ordinary graded Betti numbers of
$R/I$. For $i=0$ one has
$H_0(p,R/I)=R/I+(y_1,\dots,y_p)$ and the behavior of its  Hilbert function, that is
$\beta_{0jp}(R/I)$,  under Gr\"obner deformation has been described in  [C].  We proved in  [C, Theorem 1] that 
 $$\beta_{0jp}(R/I)\leq \beta_{0jp}(R/\ini_\tau(I))$$  for all $j$ and $p$ and for all initial ideal
$\ini_\tau(I)$ of $I$.   Note also that, since generic linear forms are an almost regular sequence  (see [AH1]), the
modules
$H_i(p,R/I)$ have
 finite length for all $i>0$ and all $p$.  We start with the following:
\medskip 

\nt {\bf Lemma 2.1} \ \  {\sl  Let $I,J$ be homogeneous ideals of $R$ and let
$\tau$ be a term order.  Then 
$$\dim_K  \Tor^R_i(R/I,R/J)_j \leq  \dim_K 
\Tor^R_i(R/\ini_\tau(I),R/\ini_\tau(J))_j
$$ for all $i$ and $j$.  } 
\medskip 

\nt The proof of Lemma 2.1 follows from a  standard  deformation argument  making use of flat families,  see [E, Chap.15]
and [Sb]   for details. We will just sketch it.  Any vector $\lambda\in \RR_+^n $ induced a graded structure on $R$ and
any monomial $m=x_1^{\alpha_1}\cdots x_n^{\alpha_n}$ is homogeneous of $\lambda$-degree $\lambda(m)=\sum \alpha_i \lambda
_i$. For every
$f\in R$ one defines its $\lambda$-degree to be the the largest degree of a monomial in $f$.    If  $f\in R$ is a
polynomial of $\lambda$-degree $a$  then one defines the initial form $\ini_\lambda(f)$ of $f$ with respect to $\lambda$
to be the sum of the terms of
$f$ of  $\lambda$-degree equal to  $a$. Let
$t$ be a new variable. If  $f=\sum \gamma_i m_i\in R$ with $\gamma_i\in K$ and $m_i$ monomials and its   
$\lambda$-degree is $a$  then one defines   the $\lambda$-homogenization $\tilde f$ of $f$ to be the polynomial 
$\tilde f=\sum \gamma_i m_i t^{a-\lambda(m_i)}$. 
Then for every ideal $I$ of $R$ one defines
$\ini_\lambda(I)$ to be the the ideal of $R$  generated by
$\ini_\lambda(f)$ with $f\in I$ and $\tilde I $ to be the the ideal of $S=R[t]$  generated by
$\tilde f$ with $f\in I$.   By construction one has $S/\tilde I \otimes S/(t)\simeq R /\ini_\lambda(I)$ and one can show
that $S/\tilde I$ is $K[t]$-free and  that $(S/\tilde I)_t\simeq R/I[t,t^{-1}]$. 

Given a term order and a finite set of ideals,
$I$ and $J$ in our case,  one can represent their initial ideals by means of a weight vector in the sense that there
exists  $\lambda\in \RR_+^n $  such that 
$\ini_\lambda(I)=\ini_\tau(I)$ an $\ini_\lambda(J)=\ini_\tau(J)$. Let $t$ be a new variable. 
Let $\tilde I  \subset S=R[t]$ and $\tilde J \subset S$ be the $\lambda$-homogenization of $I$ and $J$ with respect to 
$t$. Consider the bigraded structure on $S$ obtained by giving   degree
$(1,\lambda_i)$ to $x_i$ and degree $(0,1)$ to $t$. 
  By construction $ S/\tilde I $ and $S/\tilde J $ are bigraded $S$-modules and so is 
$T_i=\Tor^S_i(S/\tilde I,S/\tilde J)$. Let
$T_{ij}$  be the direct sum of all the components of $T_i$ of bidegree $(j,k)$ as
$k$ varies. Since
$T_{ij}$ is a finitely generated and graded $K[t]$-module we may decompose it as 

$$T_{ij}=F_{ij} \oplus G_{ij}$$   where  $F_{ij}$ is the free part  and $G_{ij}$ is the   torsion part  which, being
$K[t]$-graded, is a direct sum of modules of the form $K[t]/(t^a)$ for various
$a>0$.  Denote by $t_{ij}$, $f_{ij}$ and $g_{ij}$, respectively, the minimal number of generators of 
$T_{ij},F_{ij}$ and $G_{ij}$  as $K[t]$-modules.

Since $t$ is a regular homogeneous element over $S$, $S/\tilde I$ and $S/\tilde J$, one has a short exact sequence: 

$$0\to  \CoKer \phi_i  \to \Tor^R_i(R/\ini_\tau(I),R/\ini_\tau(J)) \to \Ker
\phi_{i-1}
 \to 0$$   where $\phi_i$ is the multiplication by $t$ from $T_i$ to itself.  It follows that the dimension of
$\Tor^R_i(R/\ini_\tau(I),R/\ini_\tau(J))$ in   degree
$j$ is given by   
$t_{ij}+g_{i-1,j}$ while the dimension of 
$\Tor^R_i(R/I,R/J)$ in degree $j$ is given by $f_{ij}$. Since
$t_{ij}=f_{ij}+g_{ij}$, we are done. \endProof
\medskip 

Let $\Delta: \RR^n \to \RR^n$ and $\delta: \RR^n \to \RR^n$  be the linear maps defined by 
$$\Delta(a)=(a_1,a_1+a_2,\dots,a_1+a_2+\dots+a_n), $$

$$\delta(a)=(a_1-a_2,a_2-a_3,\dots,a_{n-1}-a_n,a_n).$$

Fix an integer $d$ and let $M_d$ the set of the monomials of degree $d$ in $R$.  The Borel order is on $M_d$ defined as
follows. Let $x^a$ and $x^b$ be monomials in
$M_d$ with exponents
$a=(a_1,\dots,a_n)$ and $b=(b_1,\dots,b_n)$ then

$$x^a\geq _{\Borel} x^b \hbox{    iff  }  \Delta(a)\geq \Delta(b)   \hbox{  componentwise }$$   The following is a
well-known fact: 
\medskip

\nt {\bf Lemma 2.2 } \ \  {\sl 1) The Borel order is a partial order on $M_d$.  \par   
\nt 2)   Let $x^a$ and $x^b$ monomials of degree $d$.  Then $x^a>_{\Borel} x^b$ iff $x^a>_{\tau} x^b$ for all term orders
$\tau$.  } 
\medskip 

The proof of 1) is easy. To prove  2)  one notes that any term order on $M_d$ can be represented by a weight function
$w=(w_1,\dots,w_n)\in \RR^n$. But $w_i>w_{i+1}$ since we consider only term orders with $x_i>x_{i+1}$.  Then the desired
statement  follows immediately from the  following identity.   Let 
$a,w
\in\RR^n$  then 
$$a \cdot w=\Delta(a) \cdot \delta(w)$$ where $\cdot$ denotes the ordinary scalar product. 

A subset $A$ of a poset $P$ is said to be an (upper) {\it poset ideal} if for all
$a\in A$ and every $b\in P$ with $b>a$ one has $b\in A$. 

Any $g=(g_{ij})\in \GL_n(K)$  acts  on $R$ as a  $K$-algebra graded  isomorphism  by 
$$g(x_i)=\sum_{j=1}^n g_{ji}x_j.$$  An ideal $I$ is said to be  Borel-fixed if it is invariant under the action of any
upper triangular matrix. A monomial ideal $I$ is said to be strongly  stable if whenever  $x_im$ is in
$I$  for some monomial $m$ and some
$i$  then  $x_jm$  is in $I$ for all $j<i$.  Equivalently, $I$ is strongly stable iff for every
$d$  the set of the monomials of degree $d$ in $I$ form a poset ideal of $M_d$ with respect to  the Borel order.    It
turns out that the strongly stable monomial ideals are exactly the Borel-fixed ideals (in positive characteristic this is
no longer the case). 

For any homogeneous ideal $I$ and  any term order $\tau$ one can consider  the generic initial ideal
$\Gin_\tau(I)$ of $I$ with respect to $\tau$,  defined as
$\ini_\tau(g(I))$ where $g$ is a generic element in $\GL_n(K)$.  One has:
\medskip 

\nt {\bf Lemma 2.3 } \ \  {\sl Let $I$ be a homogeneous ideal and $\tau$ and
$\sigma$ be term orders. \par 
\nt 1) $\Gin_\tau(I)$ is a Borel-fixed ideal. \par 
\nt 2) $\Gin_\tau(I)=I$ iff $I$ is Borel-fixed.\par 
\nt 3) $\Gin_\sigma(\Gin_\tau(I))=\Gin_\tau(I)$. }\medskip 

\nt Proof: 1) is proved, for instance, in [E, Chap.15]. 3) follows from 1) and 2). One direction of implication 2)
follows from 1).  It remains to prove that $\Gin_\tau(I)=I$ for any  Borel-fixed ideal
$I$.   We recall that a generic matrix  $g$ has an LU decomposition, that is, it  factors  as a product $g=ab$ where
$b$ is upper triangular and $a$ is lower triangular. One can even take one of the two matrices with $1$'s  on the
diagonal but we do not mind.  This is a well-known theorem in matrix theory and a proof of it  with  a description of
the  precise   conditions that $g$ must satisfy for having such a decomposition can be found for instance in    [GV,
Theorem 3.2.1].  For every lower triangular matrix
$a$ and every monomial
$m$ one has 
$$a(m)=\lambda m+\dots\hbox{ monomials which are } <_{\Borel} m. \eqno{(1)}$$   It follows that for every polynomial 
$f$  one has 
$\ini_\tau(a(f))=\ini_\tau(f)$ and hence 
$\ini_\tau(a(J))=\ini_\tau(J)$ for every ideal $J$. Summing up, 
$$\Gin_\tau(I)=\ini_\tau(g(I))=\ini_\tau(ab(I))=\ini_\tau(a(I))=\ini_\tau(I)=I$$ where the third equality holds  since
$b(I)=I$ by assumption and the last holds since $I$ is a monomial ideal.     
\endProof 
\medskip 

\nt {\bf Lemma 2.4 } \ \  {\sl Let $I$ be a Borel-fixed ideal and let $m_1,\dots,m_k$ be its monomial generators. Let
$g\in \GL_n(K)$ be a generic matrix. Then $g(I)$ is generated by  polynomials 
$f_1,\dots,f_k$  of the form   $f_i=m_i+h_i$ such that  the monomials in $h_i$ are smaller than $m_i$  in the Borel
order.  The  polynomials $f_1,\dots,f_k$  form a   Gr\"obner basis of $g(I)$ with  respect to any term order.  }
\medskip 

\nt Proof:   As noted already in the proof of 2.3 we may write $g=ab$ with  where $b$ is upper triangular and
$a$ is lower triangular. Then $g(I)=ab(I)=a(I)=(a(m_1),\dots, a(m_k))$.  Set $f_i= \lambda_i^{-1} a(m_i)$  where
$\lambda_i$ is the leading coefficient of $a(m_i)$.   By $(1)$ we know that the $f_i$ have the desired structure and the
rest follows from Lemma 2.2.
\endProof

Proper sequences will be of importance in the next sections. Let us recall the definition from [HSV, Definition 2.1,
Remark 2.4]:

\medskip 
\nt {\bf Definition   2.5} \ \  {\sl  Let   $z_1,\dots,z_k$ be a sequence of homogeneous elements in $R$ and let
$M$ be a graded $R$-module. Set
$I=(z_1,\dots,z_k)$. We say that $z_1,\dots,z_k$ is a proper
$M$-sequence if one of the following equivalent conditions hold:\par
\nt 1) $z_{j+1}H_i(z_1,\dots,z_{j},M)=0$ for $j=1,\dots, k-1$ and $i>0$.\par
\nt  2) $IH_i(z_1,\dots,z_{j},M)=0$ for $j=1,\dots, k$ and  for $i>0$.\par } 
\medskip 

Let  $z_1,\dots,z_k$ be  a proper $M$-sequence of linear forms.  Then the long exact sequence of Koszul homologies (see
[BH, Corollary 1.6.13]) splits  into shorter ones. Denote by 
$Z_j$ the sequence $z_1,\dots, z_j$.   For all $i\geq 2$ and all $j\geq 1$ one has: 

$$0\to H_i(Z_j,M)\to H_i(Z_{j+1},M) \to H_{i-1}(Z_j,M)(-1)\to 0 \eqno{(2)}$$
 
and  
$$0\to H_1(Z_j,M)\to H_1(Z_{j+1},M) \to H_{0}(Z_j,M)(-1) \longrightarrow^{ 
\!\!\!\!\!\!\! {\bullet}}\  H_{0}(Z_j,M)\to  H_0(Z_{j+1},M) \to 0 \eqno{(3)}$$ 
  where the map 
$\longrightarrow^{  \!\!\!\!\!\!\!{\bullet}} \ $   
 is the multiplication by $z_{j+1}$ (up to sign).

\bigskip\bigskip\bigskip

\centerline{\bf 3. Koszul-Betti numbers and  Borel-fixed ideals } 
\medskip

First of all, we explain why the last $p$ variables are  ``generic"  for a Borel-fixed ideal.

\nt {\bf Lemma 3.1} \ \  {\sl  Let $I$  be a Borel-fixed ideal.  For a given $p$,
$1\leq p\leq n$, let   $X_p$ be the sequence $x_{n-p+1},x_{n-p+2},\dots,x_n$. Then 
$\beta_{ijp}(R/I)=\dim_K H_i(X_p, R/I)_j$.
  } 
\medskip 

\nt Proof:  Let $y=y_1,\dots, y_p$ be a   sequence of $p$ generic linear forms.  We can choose an upper-triangular matrix
$g$ such that the induced
$K$-algebra graded  isomorphism  $g:R\to R$ maps the linear space generated by
$y_i$ to the linear space generated by $X_p$.  Then $H_i(y,R/I)=H_i(X_p, R/g(I))$ and $g(I)=I$ by assumption.
\endProof
\medskip

For a Borel-fixed  ideal $I$ a formula for the  Koszul-Betti numbers  of $R/I$ has been described by Aramova and Herzog
[AH, Proposition 2.1].  Their result holds  more generally for  a stable ideal and can be seen as an extension of the
Eliahou-Kervaire [EK]  formula for the  Betti numbers.  To state their formula  in the form which best suits our needs we
introduce a   useful piece of notation.  For a monomial  $u$ with exponent
$a=(a_1,\dots, a_n)$ one defines 
$$\max(u)=\max\{ i : a_i>0\}.$$
  For a set of monomials $A$ and for $i=1,\dots,n$ we put:
$$m_i(A)=|\{ u \in A : \max(u)=i\} |$$
$$m_{\leq i}(A)=|\{ u \in A : \max(u)\leq i\} |$$ 
$$m_{ij}(A)=|\{ u \in A : \max(u)=i  \hbox{ and  } \deg(u)=j\} | $$  When $I$ is either a vector space generated by
monomials of the same degree or a monomial ideal,  we set  
$$m_i(I)=m_i(G), \quad m_{\leq i}(I)=m_{\leq i}(G), \quad  m_{ij}(I)=m_{ij}(G)$$   where $G$ is the set of the minimal
monomial (vector space or ideal) generators  of
$I$. 

Taking into consideration Lemma 3.1,  the result of Aramova and Herzog [AH] can be stated as follows: 

\medskip 
\nt {\bf Theorem 3.2} (Aramova-Herzog)  \ \  {\sl  Let $I$  be a  Borel-fixed ideal. Then for all
$i>0$ and for all $j$ and $p$ the Koszul-Betti numbers of $R/I$ are given by the formula: 
$$
\beta_{ijp}(R/I)=\sum_{s=i+n-p}^n m_{s,j-i+1}(I) {s+p-n-1 \choose i-1} 
$$ } 
\medskip 

With $p=n$ one gets the Eliahou-Kervaire formula for Betti numbers:  

\medskip 
\nt {\bf Theorem 3.3} (Eliahou-Kervaire)  \ \   {\sl  Let $I$  be a  Borel-fixed ideal. Then for all
$i>0$ and for all $j$  the  Betti numbers of $R/I$ are given by the formula: 
$$
\beta_{ij}(R/I)=\sum_{s=i}^n m_{s,j-i+1}(I) {s-1 \choose i-1} 
$$ } 
\medskip  

We have  also: 

\medskip 
\nt {\bf Lemma 3.4} \ \  {\sl  Let $I$  be a  Borel-fixed ideal. Then for all  
$j$ and $p$  one has: 
$$\beta_{0jp}(R/I)={n-p-1+j \choose n-p-1}-m_{\leq n-p}(I_j)$$ } 
\medskip 

\nt Proof: This  follows immediately from the definition of $m_{\leq i}$ since, by Lemma 3.1,  the sequence  
$x_{n-p+1},x_{n-p+2},\dots,x_n$  is generic for
$I$.\endProof 

The following is a useful property of Borel-fixed spaces, see [B,Proposition 1.3]. 

\medskip 
\nt {\bf Proposition  3.5} \ \  {\sl  Let $B$ be a  Borel-fixed  vector space of monomials of degree $d$. Then one has
$m_{i}(R_1B)= m_{\leq i}(B).$ } 
\medskip 

The functions $m_i$ and $m_{\leq i}$  are important for comparing Borel-fixed ideals. Here is why:

\medskip 
\nt {\bf Proposition  3.6} \ \  {\sl  Let $I,J$ be   Borel-fixed ideals with the same Hilbert function. Assume that 
$m_{\leq i}(J_j)\leq m_{\leq i}(I_j)$ for all $i$ and
$j$. Then one has: \par 
\nt  i)  $m_i(I)\leq m_i(J)$ for all $i$. \par
\nt  ii) $\beta_{ijp}(R/I)\leq \beta_{ijp}(R/J)$ for all $i,j,p$.  }\medskip 

\nt Proof:  i) Let $k$ be  an integer such that both $I$ and $J$ are generated in degree $\leq k$.  The ideals
$I_{<k>}$ and $J_{<k>}$ are Borel-fixed generated in degree $k$ and have the same Hilbert function. By the
Eliahou-Kervaire formula,  they both have  a linear resolution.  It follows that the Betti numbers of
$I_{<k>} $ and $J_{<k>}$ are the same. Again by the Eliahou-Kervaire formula it follows that
$m_i(I_k)=m_i(J_k)$ for all $i$ and hence  $m_{\leq i}(I_k)=m_{\leq i}(J_k)$ for all $i$. Now note that we have
$m_i(I)=\sum_{j=1}^k m_{ij}(I)$ and 
$m_{ij}(I)=m_i(I_j)-m_i(R_1I_{j-1})$. We may write $m_i(I_j)=m_{\leq i}(I_j)-m_{\leq i-1}(I_j)$ and by 3.5 we have 
$m_i(R_1I_{j-1})=m_{\leq i}(I_{j-1})$. Summing up, we have 
$$m_i(I)=\sum_{j=1}^k \left [ m_{\leq i}(I_j)-m_{\leq i-1}(I_j)-m_{\leq i}(I_{j-1}) \right ]= m_{\leq
i}(I_k)-\sum_{j=1}^k m_{\leq i-1}(I_j) $$   The same formula holds also when we replace
$I$ with $J$. Since we know that
$m_{\leq i}(I_k)=m_{\leq i}(J_k)$, it follows that  
$$m_i(J)-m_i(I)=\sum_{j\geq 1}  \left [ m_{\leq i-1}(I_j) -m_{\leq i-1}(J_j)
\right ]  $$    By assumption $m_{\leq i-1}(I_j) -m_{\leq i-1}(J_j)$ is non-negative and hence  we may conclude that
$m_i(J)\geq m_i(I)$. \par

ii) By virtue of Lemma 3.4,   for  $i=0$ one has 
$$\beta_{0jp}(R/J)- \beta_{0jp}(R/I)=  m_{\leq n-p}(I_j) - m_{\leq n-p}(J_j). $$   Then $\beta_{0jp}(R/I)
\leq \beta_{0jp}(R/J)$ follows from the assumption.  To prove the inequality for
$i>0$ we just follows  Bigatti's  proof of Theorem 1.1 (b)   by replacing the Eliahou-Kervaire formula    with the 
Aramova-Herzog formula.  For reader's convenience we reproduce it.  First, by 3.5, in the Aramova-Herzog formula   we may
replace  
$m_{s,j-i+1}(I)$ with $m_{s}(I_{j-i+1})-m_{\leq s}(I_{j-i})$. Then we may replace
$m_{s}(I_{j-i+1})$ with   $m_{\leq s}(I_{j-i+1})-m_{\leq s-1}(I_{j-i+1})$. We get 
$$\beta_{ijp}(R/I)=\sum_{s=i+n-p}^n  \left [ m_{\leq s}(I_{j-i+1})-m_{\leq s-1}(I_{j-i+1})-m_{\leq s}(I_{j-i})
\right ] {s+p-n-1 \choose i-1}$$  which we can rewrite as a sum of three terms: 
$$m_{\leq n}(I_{j-i+1}){p-1 \choose i-1} - m_{\leq n-p+i-1}(I_{j-i+1}) \eqno{(i)}$$
 $$
 \sum_{s=i+n-p}^{n-1}  m_{\leq s}(I_{j-i+1})\left [ {s+p-n-1 \choose i-1} -{s+p-n
\choose i-1} \right ]
\eqno{(ii)}$$
$$ -\sum_{s=i+n-p}^{n}  m_{\leq s}(I_{j-i})   {s+p-n-1 \choose i-1}.
\eqno{(iii)} $$ Note that $m_{\leq n}(I_{j-i+1})=\dim_K (I_{j-i+1})$ and hence (i) can be rewritten as
$$\dim_K (I_{j-i+1}) {p-1 \choose i-1} - m_{\leq n-p+i-1}(I_{j-i+1}). \eqno{(iv)}$$ Also, applying Pascal's triangle
formula, part (ii) is equal to 
$$-\sum_{s=i+n-p}^{n-1}  m_{\leq s}(I_{j-i+1})  {s+p-n-1 \choose i-2}.
\eqno{(v)}$$ Summing up,
$\beta_{ijp}(R/I)$ is the sum of (iii), (iv) and (v).  The crucial consequence is that
$\beta_{ijp}(R/I)$ can be written as $\dim_K (I_{j-i+1}) {p-1 \choose i-1}$ minus a linear combination with non-negative
coefficients of $m_{\leq a}(I_b)$ for various $a$ and $b$.  Since
$\dim_K (I_{j-i+1})=\dim_K (J_{j-i+1})$, it follows that 
$\beta_{ijp}(R/J)-\beta_{ijp}(R/I)$ can be written as a linear combination  with non-negative coefficients of
$m_{\leq a}(I_b)-m_{\leq a}(J_b)$. By assumption the latter is non-negative and  
 we conclude  that
$\beta_{ijp}(R/J)-\beta_{ijp}(R/I)\geq 0$. \endProof
\medskip

\medskip 
\nt {\bf Proposition  3.7} \ \  {\sl  Let $I,J$ be   Borel-fixed ideals with the same Hilbert function. Assume that 
$m_{\leq i}(J_j)\leq m_{\leq i}(I_j)$ for all $i$ and
$j$.   Then the following conditions are equivalent: \par
 i)   $\beta_{ijp}(R/I)=\beta_{ijp}(R/J)$ for all $i,j$ and  $p$,\par
 ii)  $\beta_{ij}(R/I)=\beta_{ij}(R/J)$ for all $i$ and $j$,\par
 iii)  $\beta_{1j}(R/I)=\beta_{1j}(R/J)$ for all   $j$,\par
 iv)  $\beta_{1}(R/I)=\beta_{1}(R/J)$,\par
 v) $m_{ij}(I)=m_{ij}(J)$ for all $i,j$. \par
 vi)  $m_{i}(I)=m_{i}(J)$ for all $i$. \par
 vii)  $m_i(I_j)=m_i(J_j)$ for all $i$ and $j$. \par
 viii)  $m_{\leq i}(I_j)=m_{\leq i}(J_j)$ for all $i$ and $j$.  } 
\medskip  

\medskip 

\nt Proof: The implications  $i)\Rightarrow ii)\Rightarrow iii)\Rightarrow iv)$  and 
$v)\Rightarrow vi)$ are obvious.  That  $i)\Leftrightarrow viii)$ and 
$vi)\Leftrightarrow viii)$ follows the proof of 3.6  while $iv)\Rightarrow vi)$  follows  from 3.6  since
$\beta_1(R/I)=\sum m_i(I)$.  Furthermore  $vii)
\Leftrightarrow viii)$ is easy and $vii), viii) \Rightarrow v)$ since, by 3.5,  we have
$m_{ij}(I)=m_i(I_j)-m_i(R_1I_{j-1})=m_i(I_j)-m_{\leq i}(I_{j-1})
$. 
\endProof

\bigskip\bigskip\bigskip

\centerline{\bf 4. Koszul-Betti numbers, Gin and Lex}
\medskip  

\nt In this section we prove the main results of this note.  To compare Betti numbers of Borel-fixed  and lex-segment 
ideals Bigatti   made a study of the behavior of the functions
$m_i(\dots)$  as one compares  a Borel-fixed space with its   lex-segment.  She proved in [B, Theorem 2.1] the 
following  crucial inequality.  An essentially equivalent form of the same result appears also in D.Bayer Ph.D. Thesis
[Ba, Lemma 8.1]. 

\medskip 
\nt {\bf Proposition  4.1} \ \  {\sl  Let $I$ be a  Borel-fixed   ideal  and let $L=\Lex(I)$ be the corresponding
lex-segment.  Then one has
$m_{\leq i}(L_j)\leq m_{\leq i}(I_j)$ for all $i$ and $j$. }  
\medskip 

This is the generalization of Theorem 1.1: 
\medskip

\nt {\bf Theorem 4.2} \ \  {\sl  Let $\tau$ be any term order and $I$ an homogeneous ideal, then\par 
\nt  a)   $\beta_{ijp}(R/I)\leq \beta_{ijp}(R/\Gin_\tau(I))$ for all $i,j,p$. \par
\nt b) $\beta_{ijp}(R/I)\leq \beta_{ijp}(R/\Lex(I))$ for all $i,j,p$.  }
\medskip

\nt Proof:  Let $y=y_1,\dots, y_p$ be a sequence of  generic linear forms. Let
$g\in \GL_n(K)$ such that the induced $K$-algebra graded  isomorphism
$g:R\to R$ maps $y_i$ to $x_{n-p+i}$ for $i=1,\dots,p$.  Denote by  $X_p$ the sequence
$x_{n-p+1},x_{n-p+2},\dots,x_n$. 
  Then 
$$H_i(y,R/I)=H_i(X_p, R/g(I))=\Tor^R_i(R/g(I), R/(X_p)).$$ It follows by Lemma 2.1: 
$$\beta_{ijp}(R/I)\leq \dim_K \Tor^R_i(R/\ini_\tau(g(I)), R/(X_p))_j =
\dim_K  H_i(X_p,R/\ini_\tau(g(I))_j.$$  Since the $y_i$ are generic,  $g$ can be chosen  generic as well. So
$\ini_\tau(g(I))=\Gin_\tau(I)$ and it is Borel-fixed. But then by Lemma 3.1 we know that 
$$\beta_{ijp}(R/\Gin_\tau(I))=H_i(X_p,R/\Gin_\tau(I))_j$$ and this proves a).

To  prove b), by virtue of a), we may replace $I$ with its gin, that is, we may assume that $I$ is Borel-fixed.  Then the
result follows from 3.6 and 4.1. 
\endProof
\medskip

\nt {\bf Remark 4.3} \ \  {\sl  a) Theorem 4.2 a) holds also in positive characteristic. \par  
\nt  b)   It is easy to see that $\beta_{0jp}(R/I)=\beta_{0jp}(R/\Gin(I))$ for all
$j$ and $p$, see   [C,Lemma 2]  for details. 
 } 
\medskip

Now we generalize Theorem 1.2:
\medskip 

\nt {\bf Theorem 4.4} \ \  {\sl   The following conditions are equivalent:
\par
 i)  $\beta_{ijp}(R/I)= \beta_{ijp}(R/\Gin(I))$ for all $i,j,p$.\par   ii)  $\beta_{1jn}(R/I)=
\beta_{1jn}(R/\Gin(I))$ for all $j$.\par   iii) $I$ is a componentwise linear ideal.\par   iv)  a generic sequence of
linear forms $y_1,\dots,y_n$ is a proper sequence over
$R/I$. }
\medskip

\nt Proof: i)  $\Rightarrow$  ii) is obvious  and ii)  $\Rightarrow$  iii) holds by Theorem 1.2. 
\par

We prove now that iii)  $\Rightarrow$  iv). Assuming  that $I$ is componentwise linear we have to show that 
$\mm H_i(p,R/I)=0$ for all $i>0$ and $p$ where $\mm$ denotes the homogeneous maximal ideal of $R$. First note that if
$I$ is generated in a single degree, say $d$, then so does $\Gin(I)$. Then by Theorem 3.2 we have that for
$i>0$ the homology module  $H_i(p,R/\Gin(I))$ is concentrated in a single degree, namely 
$d+i-1$. By Theorem 4.2 the same is  true also for
$H_i(p,R/ I)$ and hence $\mm H_i(p,R/ I)=0$.  Now assume that  $I$ is possibly generated in distinct   degrees.  For a
fixed $p$,   let $K_{\bullet}=K(p,R)$ the Koszul complex of $p$ generic linear forms and let
$\phi_i$ be the map from
$K_i$ to $K_{i-1}$. Let
$a\in H_i(p,R/I)$ a homogeneous element, say of degree $s$. Consider its  preimage, say, 
$a=\bar f$ with
$f\in K_i$. Then $\phi_i(f)$ is in $IK_{i-1}$. Since $\phi_i(f)$ is homogeneous of degree $s$, we deduce that
$\phi_i(f)$ is in $I_kK_{i-1}$ where $k=s-i+1$. Set 
$J=I_{<k>}$, that is $J$ is  the ideal generated by $I_k$. We have that the class of
$f$ is in $H_i(p,R/J)$. By construction $J$ is generated in a single degree and by assumption it has a linear
resolution.  From what we have seen  above we may conclude that $\mm f$ is contained in $\Image
\phi_{i+1}+ J K_i$. Since
$J\subset I$ we have that $\mm f$ is contained in $\Image \phi_{i+1}+ I K_i$ which in turns imply that $\mm a=0$ in
$H_i(p,R/I)$. 

It remains to prove that iv) implies i).   We have mentioned already in Remark 4.3 that 
$\beta_{0jp}(R/I)=\beta_{0jp}(R/\Gin(I))$ for all
$j$ and $p$ holds for any ideal $I$.  Assuming iv), it is then sufficient to show that the numbers
$\beta_{ijp}(R/I)$ only depend on  the numbers
$\beta_{0jp}(R/I)$.   For $i=1$ the exact sequence (2) implies that: 
$$\beta_{1jp}(R/I)=\beta_{1jp-1}(R/I)+\beta_{0j-1p-1}(R/I)-\beta_{0jp-1}(R/I)+
 \beta_{0jp}(R/I)$$

For $i>1$ the exact sequence (3) implies: 
$$\beta_{ijp}(R/I)=\beta_{ijp-1}(R/I)+\beta_{i-1j-1p-1}(R/I)$$

Since $\beta_{ijp}(R/I)=0$ for all $i>p$  these recursive relations  imply what  we want and  conclude the proof. 
\endProof

\medskip 
 
Now we  generalize  Theorem 1.3:

\medskip

\nt {\bf Theorem 4.5} \ \  {\sl   The following conditions are equivalent:
\par
 i)  $\beta_{ijp}(R/I)= \beta_{ijp}(R/\Lex(I))$ for all $i,j,p$.\par  ii) 
$\beta_{1jn}(R/I)=
\beta_{1jn}(R/\Lex(I))$ for all $j$.\par  iii) $I$ is a Gotzmann ideal.\par   iv) 
$\beta_{0jp}(R/I)=
\beta_{0jp}(R/\Lex(I))$ for all $j,p$ and $I$ is componentwise linear.\par  } 
\medskip 

\nt Proof: Set $L=\Lex(I)$. i)  $\Rightarrow$  ii) is obvious  and ii) 
$\Rightarrow$  iii) holds by Theorem 1.3.  A Gotzmann ideal  is clearly  componentwise linear. So to prove that iii)
implies iv), it suffices to prove  that $\beta_{0jp}(R/I)=
\beta_{0jp}(R/L)$ for all
$j,p$. By  Remark 4.3    passing to the gin does not change
$\beta_{0jp}$. So we may replace
$I$ with $\Gin(I)$ which is still Gotzmann. In other words we may assume that
$I$ is  Borel-fixed and $\beta_1(R/I)=\beta_1(R/L)$.  Then by  4.1 and  3.7 we know that
 $\beta_{0jp}(R/I)=\beta_{0jp}(R/L)$ for all $j,p$. Finally, we prove that iv) implies i).  Since $I$ is componentwise
linear, by virtue of Theorem 3.4, we may replace $I$ with its gin, that is, we may   assume that
$I$ is Borel-fixed.  The assumption
$\beta_{0jp}(R/I)= \beta_{0jp}(R/L)$, by virtue of Lemma 3.4 translates into
$m_{\leq i}(I_j)=m_{\leq i}(L_j)$ for all $i$ and $j$. By  3.7 and 4.1  this implies  
$\beta_{ijp}(R/I)= \beta_{ijp}(R/L)$ for all   $i,j,p$. 
\endProof

\nt {\bf Remark/Example  4.6} \ \  {  a) With the notation and the assumption of 3.6,  it is not true in general that
$m_{ij}(I)\leq m_{ij}(J)$ for all $i,j$. For instance take $I=(x_1,x_2)^2+(x_1,x_2,x_3)^3$.  Then $I$ is Borel-fixed, its
lex-segment  ideal is $L=x_1(x_1,x_2,x_3)+(x_1,x_2,x_3)^3$.  By 4.1 we know that $m_{\leq i}(L_j)\leq  m_{\leq i}(I_j)$
for all $i,j$ but 
$m_{22}(I)=2$ and $m_{22}(L)=1$.
\par  
\nt b) By 3.7 and 4.1 we know that a Borel-fixed ideal $I$ with lex-segment ideal $L$ is Gotzmann iff
$m_i(I)=m_i(L)$ for all $i$.  This is an example, perhaps the
 smallest one,  of a Borel-fixed Gotzmann ideal $I$  which is not a lex-segment
 ideal. In $K[x_1,x_2,x_3]$ consider 
 $$I=(x_1^3, x_1^2x_2, x_1^2x_3, x_1x_2^2, x_1x_2x_3, x_2^3, x_2^2x_3).$$ Its
 lex-segment ideal is 
$$L=(x_1^3, x_1^2x_2, x_1^2x_3, x_1x_2^2, x_1x_2x_3, x_1x_3^2, x_2^3).$$ }  To verify that $I$ is Gotzmann it suffices to
note  that  $m_1(I)=m_1(L)=1$, $m_2(I)=m_2(L)=2$ and
$m_3(I)=m_3(L)=4$. 
\medskip  

\bigskip\bigskip\bigskip 

\centerline {5. \bf  Comparing the gins  } 
\medskip 

Given a homogeneous ideal $I$ we may consider the set 
$$\GGin(I)=\{ \Gin_\tau(I) : \tau \hbox{ is a term order } \}$$ of all the generic initial ideals of $I$. Among them, the
gin-revlex (which we simply denote by $\Gin(I)$) plays a special role.  For instance,  it is known that  $\projdim(I)=
\projdim(\Gin(I))$ and $\reg(I)=\reg(\Gin(I))$,  where $\projdim$ and $\reg$ denote the projective dimension and the
Castelnuovo-Mumford regularity,  see [BS].  It follows that $\projdim(J)\geq
\projdim(\Gin(I))$ and $\reg(J)\geq \reg(\Gin(I))$ for all  
$J\in \GGin(I)$. We generalize this and show that: 

\medskip

\nt {\bf Theorem 5.1} \ \  {\sl    Let $I$ be a homogeneous ideal. Let $\tau$  be   a term order.   Then 
$$\beta_{ijp}(\Gin(I))\leq \beta_{ijp}(\Gin_\tau(I))$$ for every $i,j,p$.
 } 
\medskip 

\nt Proof:  Set $J=\Gin(I)$ and  $H=\Gin_\tau(I)$. Note  $J$ and
$H$ are Borel-fixed ideals with the same Hilbert function. By virtue of 3.6 it is therefore enough   to show that $m_{\leq
i}(H_j)\leq m_{\leq i}(J_j)$ for all $i$ and $j$. Let
$a_1\dots,a_k$ be the generators of $J_j$ and  $b_1\dots,b_k$ that of $H_j$.  We may order    the
$a_r$'s and the $b_r$'s  according to the revlex order. Since by Lemma 2.3
$\Gin(H)=H$, it follows from [C, Corollary 6] that 
$a_r\geq b_r$ in the revlex order for all $r$. This implies that $\max(a_r)\leq \max(b_r)$ for all $r$ and hence 
$m_{\leq i}(H_j)\leq m_{\leq i}(J_j)$.
\endProof 

More generally, the proof above shows that Theorem 5.1 holds also if $\Gin_\tau(I)$ is replaced by any Borel-fixed
initial ideal of $I$.  Theorem 5.1 should be compared with the example in [CE]   showing   that,  in  general,  there is
no ideal with the smallest Betti numbers in the family of ideals with a given Hilbert function and with the example in 
[F, Section 6]  showing  that there is  no ideal with the smallest Betti numbers in the family of Borel-fixed ideals with
a given  Hilbert function.    Another application of [C,Corollary 6]  yields easily: 
\medskip

\nt {\bf Proposition  5.2} \ \  {\sl  Among all the elements in $\GGin(I)$  the gin-lex is the closest to the lex-segment
ideal of $I$, in the sense that if $L$ is the lex-segment ideal of $I$, $J$ is the gin-lex of $I$ and $H$ is any other
gin of $I$  then $\dim (L\cap J)_j\geq \dim (L\cap H)_j$ for all $j$. } 
\medskip 

We know that the lex-segment ideal has the largest Betti numbers in the class of the ideals with a given Hilbert function
and that in $\GGin(I)$  the gin-lex is the closest ideal to the  lex-segment. Therefore   it makes sense to ask whether
the gin-lex has the largest Betti numbers among all the ideals in $\GGin(I)$ .  It turns out that (quite surprisingly)
this is not the case in general.  Even more interesting, there are ideals
$I$ such that there is no ideal with the  largest Betti numbers in $\GGin(I)$ and the gin-lex need not  have maximal
Betti numbers.   The main difficulty in finding  examples with these  pathologies 
 is that one must be able to detect all the gins of a given ideal which is a  hard  task.  The
examples we are going to present belong to a family of ideals called  almost Borel-fixed. 
 Let us define this family.  Let $d\in \NN$.  Given a Borel-fixed space of monomials  $A$ in $M_d$ and a  monomial  $b\in
M_d$ we say that $b$ is a lower neighbor of 
$A$ if $b\not\in A$ and $a\in A$ whenever $a>_{\Borel} b$. We denote by $\Ln(A)$ the set of the lower neighbor of $A$. 

\medskip 

\nt {\bf Definition  5.3} \ \  {Let $A$ be a Borel-fixed space of monomials in $M_d$.  Let $W$ be the vector space
generated by the elements in $\Ln(A)$ and let 
$V\subseteq  W$ be a subspace.  The vector space   $A+V$ is called an almost Borel-fixed space.  A homogeneous ideal $I$
is said to be almost Borel-fixed if for each $d\in \NN$ the space
$I_d$ is almost Borel-fixed. }  
\medskip 

The main property of  almost Borel-fixed spaces and  ideals is that one has a complete description  of the set of all the
gins.  We have: 
\medskip 
 
\nt {\bf Proposition  5.4} \ \  {\sl  Let $A$  and  $V$ be as in the Definition 5.3.    Then for every term order $\tau$
one has: 
$$\Gin_\tau(A+V)=A+\ini_\tau(V).$$ } \medskip 

\nt Proof: The left and right hand side  of the equality we have to prove are vector spaces of the same dimension. 
Therefore it is enough to prove the inclusion $\supseteq$.  Since $A+V\supseteq A$  we have  
$\Gin_\tau(A+V)\supseteq \Gin_\tau(A)=A$.  To conclude the proof we need  to show that    for a generic
$g\in \GL_n(K)$ and for  all $f \in V$ one has $\ini_\tau(f) \in \ini_\tau(g(A+V))$.  Let $m$ be a lower neighbor of $A$.
Note that $A$ and  $A\cup \{m\}$ are  Borel-fixed sets. Then, by virtue of 2.3 and 2.4,
$\ini_\tau(g(A))=A$  and the normal form of $g(m)$   with respect to $g(A)$  has the form  $m+h$ where  $h$ contains
only  monomials  which are  smaller (in the Borel order) then $m$. Now say 
$f=\lambda_1m_1+\dots+\lambda_r m_r$ where  the $m_i$ are lower neighbors of $A$ and $\lambda_i \in K$.  It follows that
the normal form of $g(f)$  with respect to $g(A)$ is  $f+H$ and each  monomial in $H$ is  smaller in the Borel order
than   some of the $m_i$. This implies that $\ini_\tau(f+H)=\ini_\tau(f)$ and hence that 
$ \ini_\tau(f) \in \ini_\tau(g(A+V))$. 
\endProof

Since $\Gin_\tau(I)=\oplus_d \Gin_\tau(I_d)$, Proposition 5.4 allows us to describe all  the gins of an almost
Borel-fixed  ideal provided, of course,  we have a description of the  decomposition of  $I_d$ as ``$A$+$V$" for each
$d$.  Note however  that if $U$ is an almost Borel-fixed space then $UR_1$ need not  be almost Borel-fixed, see the
next example.

\medskip 
\nt {\bf   Example   5.5} \ \  a) The simplest almost Borel-fixed space (which is not Borel-fixed) is the following: in
$3$ variables,  $A=\langle x_1^2, x_1x_2\rangle$, 
$\Ln(A)=\{ x_1x_3, x_2^2\}$ and $V=\langle x_1x_3+x_2^2\rangle$. Then, according to Proposition 5.4, the almost
Borel-fixed space 
$U=A+V$ has only two distinct gins, the gin-revlex $A+\langle x_2^2\rangle$  and the gin-lex  
$A+\langle x_1x_3\rangle$.\par
\nt b) if we embed  the example of part a) in a ring with an extra variable $x_4$  then 
 it  is still almost Borel-fixed but it is easy to see that $R_1U$ is not almost Borel-fixed. 
\medskip

\medskip 
\nt {\bf   Construction  5.6} \ \  {  One can construct an almost Borel-fixed ideal  as follows: Let $T$ be a set of
Borel-incomparable elements in $M_d$.  Set 
$$X =\langle  n \in M_d :  \hbox{ there exists } m\in T  \hbox{ such that  } n\geq _{\Borel} m\rangle$$
$$A =\langle  n \in M_d :  \hbox{ there exists } m\in T  \hbox{ such that  } n>_{\Borel} m \rangle. $$   Let
$f_1,\dots,f_p$  be  polynomials with disjoint supports such that each $f_i$ is a  sum of elements in $T$.   Let 
$B$ be a Borel-fixed  subspace of $M_{d+1}$ such that $B$ contains $X R_1$  (e.g. $B=XR_1$ or $B=M_{d+1}$).  Then the
ideal
$I=(A)+(f_1,\dots,f_p)+(B)$ is almost Borel-fixed and: 
$$\Gin_\tau(I)=(A)+(\ini_\tau(f_1), \dots, \ini_\tau(f_p) )+(B).$$  One can iterate  the construction by taking  a set of
incomparable elements in 
$M_{d+1}\setminus B$ and so on to get an almost Borel-fixed ideal which is non-Borel fixed  in more than one degree. 
\medskip 

\medskip 
\nt {\bf   Example 5.7} \ \  {Applying Construction 5.6 in $K[x_1,x_2,x_3,x_4]$  to $T=\{x_1x_3^2, x_2^2x_4\}$ with
$p=1$, $f=x_1x_3^2+x_2^2x_4$, and $B=XR_1$ one gets 
$$A=\langle x_1^3, x_1^2x_2, x_1^2x_3, x_1^2x_4,  x_1x_2^2, x_1x_2x_3, x_1x_2x_4, x_2^3, x_2^2x_3\rangle $$ and the ideal
$$I=(A)+(f)+(x_2^2x_4^2).$$ The ideal $I$ is almost Borel-fixed and has exactly two gins. The gin-lex (which is also the
gin-revlex in this case)    is: 
$$G_1=(A)+(x_1x_3^2, x_2^2x_4^2)$$ and the gin with respect to the term order induced, for instance,  by the weight
function $(6,5,2,1)$ is 
$$G_2=(A)+(x_2^2x_4, x_1x_3^3, x_1x_3^2x_4 ).$$ The Macaulay diagrams of the Betti numbers of $G_1$ and $G_2$ are
respectively 
$$\matrix{
  10 & 17 & 10 & 2 \cr 
  1  &  3 & 3  & 1} \qquad \qquad    
 \matrix{ 
  10& 18& 12& 3 \cr
  2& 5& 4& 1  } $$ This shows that the gin-lex need not  have the largest Betti numbers among all the gins of a given
ideal.  } 

In order to describe an ideal $I$ such that in $\GGin(I)$ there is no  ideal with largest Betti  numbers we need to
enlarge a little the ambient space, i.e. we need  more variables and higher degree monomials.  
\medskip 

\nt {\bf   Example 5.8} \ \  {Applying Construction 5.6 in $K[x_1,\dots,x_7]$  to $$T=\{x_1x_3x_6^2, x_2^2x_3x_7,
x_1x_4^2x_6, x_2^2x_4^2\}$$    with $p=2$, 
$f_1=x_1x_3x_6^2+x_2^2x_3x_7$ and $f_2=x_1x_4^2x_6+x_2^2x_4^2$, and
$B=XR_1$ one gets 

$A=\langle x_1^4,x_1^3x_2,x_1^3x_3,x_1^3x_4,x_1^3x_5,x_1^3x_6,x_1^3x_7,x_1^2x_2^2,x_1^2x_2x_3,
x_1^2x_2x_4,x_1^2x_2x_5,x_1^2x_2x_6,$

$x_1^2x_2x_7,x_1^2x_3^2,x_1^2x_3x_4,x_1^2x_3x_5,x_1^2x_3x_6,x_1^2x_3x_7,x_1^2x_4^2,x_1^2x_4x_5,  x_1^2x_4x_6,x_1^2x_5^2,$ 

$x_1^2x_5x_6, x_1^2x_6^2,x_1x_2^3,x_1x_2^2x_3,x_1x_2^2x_4,x_1x_2^2x_5,x_1x_2^2x_6,x_1x_2^2x_7, x_1x_2x_3^2, x_1x_2x_3x_4,$

$x_1x_2x_3x_5, x_1x_2x_3x_6, x_1x_2x_3x_7,x_1x_2x_4^2, x_1x_2x_4x_5,x_1x_2x_4x_6,x_1x_2x_5^2, x_1x_2x_5x_6,$

$x_1x_2x_6^2,x_1x_3^3,x_1x_3^2x_4,x_1x_3^2x_5,x_1x_3^2x_6,x_1x_3x_4^2,x_1x_3x_4x_5,x_1x_3x_4x_6,  x_1x_3x_5^2,$

$x_1x_3x_5x_6,x_1x_4^3,x_1x_4^2x_5,x_2^4,x_2^3x_3,x_2^3x_4,x_2^3x_5,
 x_2^3x_6,x_2^3x_7,x_2^2x_3^2,x_2^2x_3x_4,x_2^2x_3x_5, $

$x_2^2x_3x_6 \rangle$

and  
$$I=(A)+(f_1,f_2)+(x_1x_4^2x_6x_7, x_2^2x_3x_7^2, x_1x_4^2x_6^2).$$

The ideal $I$ is almost Borel-fixed and has exactly $3$ gins. The gin-revlex  
$$G_1=(A)+(x_2^2x_4^2,x_1x_3x_6^2)+( x_1x_4^2x_6x_7, x_2^2x_3x_7^2, x_1x_4^2x_6^2),$$ the gin-lex
$$G_2=(A)+(x_1x_4^2x_6, x_1x_3x_6^2)+( x_2^2x_4^2x_7, x_2^2x_3x_7^2,  x_2^2x_4^2x_6, x_2^2x_4^2x_5, x_2^2x_4^3)$$ and the
gin w.r.t. to the term order induced by the weight $(7,6,5,4,3,2,1)$
$$G_3=(A)+(x_2^2x_4^2, x_2^2x_3x_7)+( x_1x_4^2x_6x_7, x_1x_3x_6^2x_7, x_1x_4^2x_6^2, x_1x_3x_6^3).$$

Here the crucial observation is that there are no term orders such the 
$\ini(f_1)=x_2^2x_3x_7$  and  $\ini(f_2)=x_1x_4^2x_6$.  The Macaulay diagrams of the Betti numbers of
$G_1,G_2$ and $G_3$  are respectively:  

$$ 
\matrix{
  64& 240& 397& 363& 190& 53& 6\cr
  3&   17& 40&  50&  35&  13& 2 } 
$$

$$
\matrix{
  64& 242& 404& 372& 195& 54& 6\cr
   5& 24& 49& 55& 36& 13& 2 }$$

$$
\matrix{
  64& 241& 402& 373& 200& 58& 7\cr
   4& 22&  50&  60&  40&  14& 2 }
$$  Therefore there is no gin with the largest Betti numbers since $G_2$ and $G_3$ have
 maximal and incomparable Betti diagrams.  }
\medskip 
 
A slight variation in the construction above allows also to describe ideals such that the gin-lex  does not even have the
largest Castelnuovo-Mumford regularity  among all the gins.  One may ask whether ideals with few generators can  behave
as the ideals of  Examples 5.7 and 5.8.  There are experimental evidences  that even a monomial ideal with $2$ generators
can have these behaviors.  For instance, we have  run CoCoA to get   the   gins  of the ideal $I=(x_1^4,x_2^2x_3x_4)$. 
The computations show  that the gin-lex  of $I$  does not have  maximal Betti numbers and that the  there seems to
be $351$ gins with $3$ different set of maximal Betti numbers.  The problem with these  computations is twofold:  first
when we compute $\ini_\tau(g(I))$ by taking a random (upper triangular) matrix  $g$ 
 what we get is  $\Gin_\tau(I)$  but only ``with high probability", and, second,   it is not clear how can we be sure that
the list of the gins is complete.  
 
Let us conclude the section by showing that: 
\medskip 

\nt {\bf Proposition  5.9} \ \  {\sl An almost Borel-fixed ideal is componentwise linear.}  
\medskip 

\nt Proof:  Let $\tau$ be the revlex order. Let $A$, $V$ be as in Definition 5.3.  Let $f_1,\dots, f_p$ be  generators of
$V$. It is harmless to    assume that the monomials $\ini_\tau(f_i)$ are distinct. Let $D$ be the set of the monomials
in $A$ and set
$G=D\cup \{f_1,\dots, f_p\}$.  Let $I$ be the  ideal  generated by the set $G$. We have to show that $I$ has a linear
resolution. Set $m_i=\ini_\tau(f_i)$ and let $J$ be the ideal generated by 
$D$ and by $m_1,\dots,m_p$. Note that $J$ has a linear resolution since  it is a  Borel-fixed ideal generated in a single
degree. Therefore it suffices to show that  $\ini_\tau(I)=J$ or, equivalently, that $G$ is a Gr\"obner basis. We apply the
Buchberger algorithm to $G$. Set  $\ini_\tau(G)=D\cup \{m_1,\dots, m_p\}$. It is suffices to show that any $S$ -polynomial
which correspond to a minimal syzygy among  elements of $\ini_\tau(G)$ reduces to $0$ via $G$.  The elements of 
$\ini_\tau(G)$ generate  a Borel-fixed space. Their minimal syzygies are described by the Eliahou-Kervaire resolution. 
They have the following form: For every term $n\in
\ini_\tau(G)$  and any $j<i=\max(n)$ set $n_1=nx_j/x_i$.  Then 
$n_1\in \ini_\tau(G)$  and $x_jn-x_in_1=0$ is a  minimal syzygy.  If $n\in D$ then $n_1\in D$ as well  and the
corresponding $S$-polynomial is $0$  as any 
$S$-polynomial among monomials. If instead $n\not \in D$, then $n=m_k$ for some $k$, say $k=1$.  In this case
$n_1$ must be in $D$ since by construction $n_1>_{\Borel} n$ and, by the very definition of almost Borel-fixed space, $n$
is a lower neighbor of $A$. So the $S$-polynomial corresponding to the given syzygy is
$S=x_jf_1-x_in_1$. Let  $b$ be a non-initial monomial of $f_1$ and set $u=\max(b)$.   Then 
$i=\max(m_1)\leq \max(b)=u$ (here we use the fact that the term order is revlex) and the 
$x_jb=x_u b_1$  with $b_1=x_j a/x_u \in D$ since $b$ is also a lower neighbor of $A$ and 
$b_1>_{\Borel} b$. This implies that $S$, being  a sum of  multiples of  monomials of $D$,   reduces to $0$ via  $G$.  
\endProof 

Note that, with the notation of the proof of 5.9, we did not prove that $\Gin(I)=J$, just that 
$\ini(I)=J$.  But since, a posteriori, we know that $I$ is componentwise linear then its gin-revlex must be generated in
the initial degree and then  it follows from  Proposition 5.4 that $\Gin(I)=J$.  Simple examples show that that, in
general, the  non-revlex  gins of $I$ might  need generators in higher degrees.

\bigskip\bigskip\bigskip 

\centerline {6. \bf  Some Questions  } 
\medskip 

We have shown that, passing to a generic initial ideal,  the Koszul-Betti numbers can only increase.  Most likely the
same is true  for any initial ideal. 
\medskip 

\nt {\bf   Question  6.1} \ \  {Is it true that $\beta_{ijp}(R/I)\leq \beta_{ijp}(R/\ini_\tau(I))$ for all
$i,j,p$ and for all term orders $\tau$?  } 
\medskip 

How much of what we have shown can be extended to positive characteristic? For instance:
\medskip   

\nt {\bf   Question  6.2} \ \  {Do Theorem 4.2(b) and Theorem 5.1 hold if the base field has positive characteristic?   } 
\medskip 

We have seen that the gin-lex need  not  have the largest Betti numbers among all the gins of a given  ideal. But, of
course, it does it whenever  it is the lex-segment. On the other hand, the gin-lex  is very rarely equal to the
lex-segment, even for ideals generated by generic polynomials. For instance the gin-lex of  two generic quadrics in four
variables differs from the lex-segment already in degree $4$.  However there are  experimental evidences that the answer
to the following questions might be positive: 
\medskip   

\nt {\bf   Question  6.3} \ \   {Let $I$ be the  ideal generated by  $m$ generic forms of degree $d$ in $n$ variables
with $m\geq n$. 
 Is it true the gin-lex of $I$ equals the lex-segment?   } 
\medskip   

\nt {\bf   Question  6.4} \ \   {Let $I$ be an ideal generated by generic forms of given degrees. Does the gin-lex have
the largest Betti numbers among all the gins of $I$?} 
\medskip 

Note that the  set $\GGin(I)$ can be very large even for ``small" ideals. For instance in four variables, for    two
generic cubics  we have detected   $93$ distinct gins (and most likely there are no others) while for two generic
quartics    we have detected more than $3000$ gins (and most likely there are many others).
In this case we have checked that the gin-lex have the largest Betti numbers among all the gins  we have found.  

In view of Theorem 1.2 and Theorem 1.3, one could ask what  happens if we assume that 
$\beta_i(R/I)=\beta_i(R/\Gin(I))$ or $\beta_i(R/I)=\beta_i(R/\Lex(I))$ for some $i>0$.  We cannot conclude that all the
Betti numbers are equal.   This is because some Betti numbers (typically  the last)  can be forced by the Hilbert
function.  For instance all the rings  with Hilbert function $1,3,4,0$   have the same last Betti number but the rest of
the resolution can vary.  On the other hand there might be some sort of rigidity toward  the end of the resolution and
the following might  be true. 
\medskip 

\nt {\bf   Question  6.5} \ \   { Let $I$ be an homogeneous ideal and $J$ be either $\Gin(I)$ or $\Lex(I)$.  Assume  that
$\beta_i(R/I)=\beta_i(R/J)$ for some $i>0$.  Does this imply that 
$\beta_j(R/I)=\beta_j(R/J)$  for all $j\geq i$?  }  
\medskip 

If we assume that $I$ is Borel-fixed and $J$ is $\Lex(I)$ then it follows easily from the  proof of 3.6 that
$\beta_i(R/I)=\beta_i(R/J)$ implies $m_j(I)=m_j(J)$ for all $j\geq i$. But then
 $\beta_j(R/I)=\beta_j(R/J)$ for all $j\geq i$ follows  from  the Eliahou-Kervaire formula.   
 
\bigskip

\nt{\bf Thanks:}\   {  The results, the examples  and the open questions  presented in  this 
the paper  have been inspired and suggested by  computations   performed  by   the
computer algebra system CoCoA  [CNR].
}

\bigskip

{\baselineskip=12pt
\vglue .8cm \centerline {\bf References} \vglue .2cm  \frenchspacing

\item{[AH]} A.Aramova, J.Herzog,    {\sl  Koszul cycles and Eliahou-Kervaire type resolutions},  J. Algebra 181 (1996),
no. 2, 347--370. 

\item{[AH1]} A.Aramova, J.Herzog,    {\sl Almost regular sequences and Betti numbers   },  Amer. J. Math. 122
    (2000), no. 4, 689--719.   

\item{[AHH]} A.Aramova, J.Herzog, T.Hibi,  {\sl  Ideals with stable Betti numbers},  Adv. Math. 152 (2000), no. 1, 72--77.

\item{[B]} A.Bigatti, {\sl Upper bounds for the Betti numbers of a given Hilbert function},  Comm. Algebra 21 (1993), no.
7, 2317--2334.

\item{[Ba]} D.Bayer, {\sl The division algorithm and the Hilbert scheme},  Ph.D. Thesis, Harvard University, June 1982,
168 pages.   

\item{[BS]} D.Bayer, M.Stillman, {\sl A criterion for detecting $m$-regularity},  Invent. Math. 87 (1987), no. 1, 1--11.  

\item{[BH]} W.Bruns, J.Herzog, {\sl  Cohen-Macaulay rings}, 
 Cambridge Studies in Advanced Mathematics, 39. Cambridge University Press, Cambridge, 1993.

\item{[CNR]} A.Capani,  G.Niesi,  L.Robbiano,  {\sl  CoCoA, a system for doing  Computations in Commutative Algebra}, 
Available   ftp from {\tt  cocoa.dima.unige.it}. 

\item{[CE]} H.Charalambous,  G.Evans, {\sl    Resolutions with a given Hilbert function},  Commutative algebra: syzygies,
multiplicities, and birational algebra (South Hadley, MA, 1992), 19--26,  Contemp. Math., 159,  Amer. Math. Soc.,
Providence, RI, 1994. 

\item{[C]} A.Conca,  {\sl  Reduction numbers and initial ideals}, preprint 2001,  to appear in the  Proc. Amer. Math. Soc.

\item{[E]} D.Eisenbud, {\sl Commutative algebra. With a view toward algebraic geometry},  Graduate Texts in Mathematics,
150. Springer-Verlag, New York, 1995. 

\item{[EK]}  S.Eliahou, M.Kervaire,  {\sl  Minimal resolutions of some monomial ideals},  J. Algebra 129 (1990), no. 1,
1--25. 

\item{[F]}    C.Francisco,   {\sl Minimal graded Betti numbers and stable ideals}, preprint 2002.   
 
\item{[GV]}    G.Golub,  C.Van Loan,  {\sl Matrix computations},  Second edition. Johns Hopkins Series in the
Mathematical Sciences, 3. Johns Hopkins University Press, Baltimore, MD, 1989. xxii+642 pp.

\item{[G]} M.Green,  {\sl Generic initial ideals},  Six lectures on commutative algebra  (Bellaterra, 1996), 119--186,
Progr. Math., 166, Birkh\"auser, Basel, 1998. 

\item{[GS]} M.Green, M.Stillman,  {\sl A tutorial on generic initial ideals },  Gr\"obner bases and applications (Linz,
1998), 90--108, London Math. Soc. Lecture Note Ser.,  251, Cambridge Univ. Press, Cambridge, 1998.    

\item{[HH]} J.Herzog, T.Hibi,  {\sl Componentwise linear ideals},  Nagoya Math. J. 153 (1999), 141--153.

\item{[HSV]} J.Herzog, A.Simis,   W.Vasconcelos,    {\sl Approximation complexes of blowing-up rings. II},  J. Algebra 82
(1983), no. 1, 53--83.

\item{[H]} H.Hulett,  {\sl Maximum Betti numbers of homogeneous ideals with a given Hilbert function},   Comm. Algebra 21
(1993), no. 7, 2335--2350.

\item{[KR}] M.Kreuzer, L.Robbiano, {\sl Computational commutative algebra 1},   Springer-Verlag, Berlin, 2000. x+321 pp.  

\item{[P1]} K.Pardue, {\sl Non-standard Borel-fixed ideals}, Ph.D. Thesis, Brandeis  University, May 1994,   69 pages.

\item{[P2]} K.Pardue, {\sl Deformation classes of graded modules and maximal Betti numbers},  Illinois J. Math. 40
(1996), no. 4, 564--585.

\item{[Sb]} E.Sbarra, {\sl  Upper bounds for local cohomology for rings with given Hilbert function},   Comm. Algebra 29
(2001), no. 12, 5383--5409. 

\item{[S]} B.Sturmfels, {\sl  Gr\"obner bases and convex polytopes }, University Lecture Series, 8. American Mathematical
Society, Providence, RI, 1996. xii+162 pp. 

\item{[V]}    G.Valla,  {\sl  Problems and results on Hilbert functions of graded algebras}. Six lectures on commutative
algebra (Bellaterra, 1996), 293--344, Progr. Math., 166, Birkh\"auser, Basel,  1998.   
 
}

\end